\documentclass[12pt]{article}
\usepackage{amssymb,amsmath,bm, url}
\usepackage{amsfonts,amsmath,latexsym,amssymb,mathrsfs, color}

\usepackage{mathrsfs}
\newcommand{\msk}{\medskip}

\begin{document}
\newcommand{\bea}{\begin{eqnarray}}
\newcommand{\ena}{\end{eqnarray}}
\newcommand{\beas}{\begin{eqnarray*}}
\newcommand{\enas}{\end{eqnarray*}}
\newcommand{\beq}{\begin{equation}}
\newcommand{\enq}{\end{equation}}
\newcommand{\bml}{\begin{multline}}
\newcommand{\eml}{\end{multline}}
\newcommand{\bmls}{\begin{multline*}}
\newcommand{\emls}{\end{multline*}}
\def\qed{\hfill \mbox{\rule{0.5em}{0.5em}}}
\newcommand{\bbox}{\hfill $\Box$}
\newcommand{\ignore}[1]{}
\newcommand{\wtilde}[1]{\widetilde{#1}}
\newcommand{\qmq}[1]{\quad\mbox{#1}\quad}
\newcommand{\qm}[1]{\quad\mbox{#1}}
\newcommand{\nn}{\nonumber}
\newcommand{\Bvert}{\left\vert\vphantom{\frac{1}{1}}\right.}
\newcommand{\To}{\rightarrow}
\def\e{{\mathbb E}}
\def\p{{\mathbb P}}

\newtheorem{theorem}{Theorem}[section]
\newtheorem{corollary}{Corollary}[section]
\newtheorem{conjecture}{Conjecture}[section]
\newtheorem{proposition}{Proposition}[section]
\newtheorem{lemma}{Lemma}[section]
\newtheorem{definition}{Definition}[section]
\newtheorem{example}{Example}[section]
\newtheorem{remark}{Remark}[section]
\newtheorem{case}{Case}[section]
\newtheorem{condition}{Condition}[section]
\newcommand{\pf}{\noindent {\it Proof:} }
\newcommand{\proof}{\noindent {\it Proof:} }

\newcommand{\RR}{\mathbb{R}}
\newcommand{\NN}{\mathbb{N}}
\newcommand{\II}{\mathbb{I}}

\title{Stein's method for the Beta distribution and the P\'olya-Eggenberger Urn}
\author{Larry Goldstein\thanks{University of Southern California, research partially supported by NSA
grant H98230-11-1-0162.} and Gesine Reinert\thanks{Oxford University}}

\maketitle

\begin{abstract}
Using a characterizing equation for the Beta distribution, Stein's method is applied to obtain bounds of the optimal order for the Wasserstein distance between the distribution of the scaled number of white balls drawn from a P\'olya-Eggenberger urn and its limiting Beta distribution. The bound is computed by making a direct comparison between characterizing operators of the target and the Beta distribution, the former derived by extending Stein's density approach to discrete distributions.

In addition, refinements are given to D\"obler's result \cite{doblerArcSine} for the Arcsine approximation for the fraction of time a simple random walk of even length spends positive, and so also to the distributions of its last return time to zero and its first visit to its terminal point, by supplying explicit constants to the present Wasserstein bound and also demonstrating that its rate is of the optimal order.
\end{abstract}

\section{Introduction}

The classical P\'olya-Eggenberger urn at time zero contains $\alpha \ge 1$
white and $\beta \ge 1$ black balls, and at every positive integer time a ball is chosen uniformly from the urn, independently of the past, and replaced along with $m \ge 1$ additional balls of the same color. With ${\cal L}$ indicating distribution, or law, and $\rightarrow_d$ indicating convergence in distribution, it is well known, see \cite{mahm} for instance, that if $S_n=S_n^{\alpha,\beta,m}$ is the number of white balls drawn from the urn by time $n=0,1,2,\ldots$ then as $n \rightarrow \infty$
\bea \label{def:Wn}
{\cal L}(W_n) \rightarrow_d {\cal B}(\alpha/m,\beta/m) \qmq{where} W_n=\frac{S_n}{n}.
\ena
Here, for positive real numbers $\alpha$ and $\beta$ we let ${\cal B}(\alpha,\beta)$ denote the Beta distribution
having density
\bea \label{Beta.density}
p(x; \alpha, \beta) &=&  \frac{x^{\alpha -1 } (1-x)^{\beta - 1}}{B(\alpha, \beta)}{\bf 1}_{\{x \in [0,1]\}},
\ena
where $B(\alpha,\beta)=\Gamma(\alpha)\Gamma(\beta)/\Gamma(\alpha+\beta)$ is the Beta function as expressed in terms of the Gamma function $\Gamma(x)$.

Using Stein's method we derive an order $O(1/n)$ bound in the Wasserstein distance $d_W$, defined in (\ref{def:Wass}), between $W_n$ and its limiting Beta distribution in (\ref{def:Wn}). We show in Remark \ref{optimal} that the rate of Theorem \ref{urntobeta} cannot be improved. Let $x \wedge y$ and $x \vee y$ denote the minimum and maximum of two real numbers $x$ and $y$, respectively.

\begin{theorem} \label{urntobeta}
For $\alpha \ge 1$ and $\beta \ge 1$ let $S_n$ be the number of white balls in $n$ draws from a P\'olya-Eggenberger urn that initially contains $\alpha$ white and $\beta$ black balls. Then with $W_n=S_n/n$ and $Z \sim {\cal B}(\alpha/m, \beta/m)$,
\beas
d_W(W_n, Z)   &\le&
\left( \frac{  m + \alpha \vee \beta}{2nm} + \frac{\alpha \beta}{nm(\alpha+\beta)} \right) (b_0 + b_1)  + \frac{3}{2n} ,
\enas
where $b_0=b_0(\alpha/m, \beta/m)$ and $b_1=b_1(\alpha/m, \beta/m)$ are given in Lemma \ref{lem:allbounds}.
\end{theorem}

Connections between Theorem \ref{urntobeta} and the work \cite{dobler} of D\"obler are spelled out in Remark \ref{dobler}.

The ${\cal B}(1/2,1/2)$ distribution, also known as the Arcsine law, describes the asymptotic distribution of many quantities that arise naturally in the study of the simple symmetric random walk $T_n=X_1+\cdots X_n$, where $X_1,\ldots,X_n$ are independent
 variables taking the values 1 and $-1$ with probability $1/2$. For instance, let $L_{2n}$ be the random variable
\beas 
L_{2n}=\sup\{m \le 2n: T_m=0 \}
\enas
giving the last return time to zero up to time $2n$. Then, see \cite{fellerI},
\bea \label{Ln:dist}
P(L_{2n}=2k)=u_{2k}u_{2n-2k} \qmq{where} u_{2m}=P(T_{2m}=0),
\ena
where $P(T_{2m}=0)=2^{-2m}{2m \choose m}$, the probability that the walk returns to zero at time $2m$.

In the limit, $(2n)^{-1}L_{2n}\rightarrow Z$ in probability,
where $Z$ has the Arcsine distribution. It is often noted that this limiting result is somewhat counter intuitive in that the Arcsine density has greatest mass near the endpoints, and least mass in the center of the unit interval, whereas in a fair coin tossing game one might assume that players are more likely to spend equal time in the lead. Perhaps at least as remarkable is the fact that the number $U_{2n}$ of segments of the walk that lie above the $x$ axis, and $R_{2n}$, the first time the walk visits the terminal point $S_{2n}$, are all equal in distribution to ${\cal L}(L_{2n})$; see \cite{fellerI} for a nice exposition.

Using the methods presented here that were available in a preprint of this article, D\"obler \cite{doblerArcSine} presents a Wasserstein bound of order $O(1/n)$ without explicit constants between the distribution of $(2n)^{-1}U_{2n}$ and the limiting Arcsine. In Section \ref{arc:sec}, essentially by applying the bounds in Lemma \ref{lem:allbounds}, we are able to attach concrete constants to the result of \cite{doblerArcSine}, as well as show the rate of the bound is optimal.
\begin{theorem} \label{arc:thm}
Let $L_{2n}$ be the last return time to zero of a simple symmetric random walk of length of length $2n$ and let $Z$ have the Arcsine distribution. Then
\beas
d_W\left(\frac{L_{2n}}{2n},Z\right) \le \frac{27}{2n}+\frac{8}{n^2}.
\enas
The same bound holds with $L_{2n}$ replaced by $U_{2n}$ or $R_{2n}$. The $O(1/n)$ rate of the bound cannot be improved.
\end{theorem}

Beginning with the introduction by Stein \cite{Stein72} of a `characterizing equation' method for developing bounds in normal approximation, to date the method has been successfully applied to a large number of the classical distributions, including the Poisson \cite{agg}, \cite{bhj}, Multinomial \cite{Loh92}, Gamma \cite{luk94},\cite{NouPec}, Geometric \cite{Pe96}, Negative Binomial \cite{BrPh99} and Exponential \cite{ChFuRo11}, \cite{PeRo11a}, \cite{PeRo11b}, as well as to non classical distributions such as the PRR family of \cite{prr} also based on P\'olya type urn models. Here we further extend the range of Stein's method by including the Beta distribution, focusing on its role as the limiting law of the fraction of white balls drawn from the P\'olya-Eggenberger urn.

The application of Stein's method here differs from the way it is usually applied in that we focus on the approximation of particular distributions whose exact forms are known, rather than develop a bound that applies to a class of complex distributions obtained by, say, summing random variables that obey weak moment and dependence conditions, as in the case of the central limit theorem. And indeed, though explicit formulas exist for the distributions we study, the need for their approximation arises regardless, as is the case also for, say, the ubiquitous use of the normal approximation for the binomial.

Urn models of the classical type, and generalizations including drawing multiple balls or starting
new urns, have received considerable attention recently; see for example \cite{argientoetal},  \cite{chauvinetal} and  \cite{chunghandjanijungreis}. Interest has partly been sparked by the ability of urn models to exhibit power-law limiting behaviour, which in turn has been a focus of network analysis, see for example \cite{durrettbook} and \cite{prr}. Connections between urn models and binary search trees are clearly explained in \cite{mahmoudsurvey}. In particular, let $m=1$ and consider the initial state of the P\'olya-Eggenberger urn
as a rooted binary tree having $\alpha$ white and $\beta$ black leaves, or external nodes.
At every time step one external node is chosen, uniformly, to duplicate, yielding a pair of leaves of the same colour. That is, the chosen external node becomes an internal node while two external nodes of the chosen colour are added.
The rule for adding an additional white leaf to the tree at time $n$ clearly is the same as the rule for adding an additional white ball to the P\'olya-Eggenberger urn for the case $m=1$, and hence the number of white leaves of the tree and white balls in the urn have same distribution. Many variations and extensions on this theme are possible. Another line of interest comes from edge reinforced random walks, because   an infinite system of independent  P\'olya-Eggenberger urns can be used to represent edge reinforced random walks on trees, see \cite{pemantlesurvey}.

\section{Characterizing equations and generators}
Stein's method for distributional approximation is based on a characterization of the target approximating distribution. For the seminal normal case considered in \cite{Stein72}, it was shown that a variable $Z$ has the standard normal distribution if and only if
\bea \label{char.normal}
\e [Zf(Z)]=\e [f'(Z)]
\ena
for all absolutely continuous functions $f$ for which $\e | f'(Z)| < \infty$.
If a variable $W$ has an approximate normal distribution, then one expects that it satisfies (\ref{char.normal}) approximately. More specifically, if one wishes to test the difference between the distribution of $W$ and the standard normal $Z$ on a function $h$, then instead of computing $\e h(W)-Nh$, where $Nh=\e h(Z)$, one may set up a `Stein equation'
\bea \label{eq.normal}
f'(w)-wf(w)=h(w)-Nh
\ena
for the given $h$, solve for $f(w)$, and, upon replacing $w$ by $W$ in (\ref{eq.normal}), calculate the expectation of the right hand side by taking expectation on the left. At first glance it may seem that doing so does not make the given problem any less difficult. However, a number of techniques may be brought to bear on the quantity $\e[f'(W)-Wf(W)]$. In particular, this expression contains only the single random variable $W$, in contrast to the difference of the expectations of $h(W)$ and $h(Z)$, depending on two distributions.

To obtain our result, we compute the distance between the distribution of the fraction of white balls drawn from the P\'olya-Eggenberger Urn and the Beta by comparing the operators that characterize them. Our approach in characterizing the urn distribution stems from what is known as the density method; see for instance, \cite{Stein2004}, \cite{ReinertSinganotes} or Section 13.1 of \cite{cgs}.
In particular, recognizing the
$-w$ in (\ref{eq.normal})
 as the ratio of $\phi'(w)/\phi(w)$ where $\phi(w)$ is the standard normal density, one hopes to replace the term $-w$ by the ratio $p'(w)/p(w)$ when developing the Stein equation to handle the distribution with density $p(w)$, and to apply similar reasoning when the distribution under study is discrete. Use of the density method in the discrete case, followed by the application of a judiciously chosen transformation, leads to the characterization of the P\'olya-Eggenberger Urn distribution given in Lemma \ref{lem:stein:urn}.

Another approach to construct characterizing equations is known as the generator method. A number of years following the publication of \cite{Stein72}, the relationship between the characterizing equation (\ref{char.normal}) and the generator of the Ornstein-Uhlenbeck process
\beas
{\cal A}f(w)=f''(w)-wf'(w),
\enas
of which the normal is the unique stationary measure, was recognized in \cite{Barbour1990}, where
it was noted that
that in some generality the process semi-group may be used to solve the Stein equation (\ref{eq.normal}).
Given this connection between Stein characterizations and generators it is natural to consider a stochastic process which has the given target as its stationary distribution when extending Stein's method to handle a new distribution.

Regarding the use of this `generator' method for extending the scope of Stein's method to the Beta distribution, we recall that the Fisher Wright model from genetics, originating in the work in \cite{Fisher1930}, \cite{Wright1945} and \cite{Wright1949}, is a stochastic process used to model genetic drift in a population and has generator given by
\beas 
{\cal A}f(x)=w(1-w)f''(w) + (\alpha(1-w) -\beta w)f'(w)
\enas
for positive $\alpha$ and $\beta$, and that the ${\cal B}(\alpha,\beta)$ distribution is its unique stationary distribution. In particular, with $Z \sim {\cal B}(\alpha,\beta)$ we have $\e {\cal A}f(Z)=0$. Let ${\cal B}_{\alpha,\beta}h=\e h(Z)$, the ${\cal B}(\alpha,\beta)$ expectation of a function $h$; we drop the subscripts when the role of the parameters $\alpha$ and $\beta$ is clear. As $\e h(Z)-{\cal B}h$ is also zero,
we are led to consider a Stein equation for the Beta distribution of the form
\bea \label{stein:beta}
w(1-w)f'(w) + \left(  \alpha(1-w) - \beta w \right) f(w) = h(w)-{\cal B}_{\alpha,\beta}h.
\ena

Lemma \ref{lem:stein:urn} provides a characterizing equation for the P\'olya urn distribution that is parallel to equation (\ref{stein:beta}). Taking differences then allows us to estimate the expectation of the right hand side of (\ref{stein:beta}) when $w$ is replaced by $W_n$ by exploiting the similarity of the two
 characterizing operators; a similar argument can be found in \cite{eichelsbacherreinert} and  \cite{holmes} for stationary distributions of birth-death chains. The results most closely related to the present work is \cite{dobler}, and its connections to the present manuscript are discussed in Remark \ref{dobler}

First we introduce some notation. We say a subset $I$ of the integers $\mathbb{Z}$ is a finite integer interval if $I=[a,b] \cap \mathbb{Z}$ for $a,b \in \mathbb{Z}$ with
$a \le b$, and an infinite integer interval if either
 $I=(-\infty,b] \cap \mathbb{Z}$ or $I=[a, \infty)  \cap \mathbb{Z}$ or $I= \mathbb{Z}$. For a real valued function $f$ let $\Delta f(k) = f(k+1) - f(k)$, the forward difference operator, and
for a real valued function $p$ taking non-zero values in the integer interval $I$ let
\bea \label{def.psi.discrete}
\psi (k)=
\Delta p(k)/p(k) \qmq{for} k \in I.
\ena

For $Z$ a random variable having probability mass function $p$ with support an integer interval $I$, let ${\cal{F}}(p)$ denote the set of all real-valued functions $f$ such that either $\e \Delta f(Z-1) $ or
$\e \psi(Z)f(Z) $ is finite,
$\lim_{n \rightarrow \infty} f(n)p(n+1) =0$ when $\sup\{k:p(k)>0\}=\infty$, $\lim_{n \rightarrow -\infty} f(n)p(n+1) =0$ when $\inf\{k:p(k)>0\}=-\infty$,
and $f(a-1)=0$ in the case where $I=[a,b] \cap \mathbb{Z}$ or $I=[a,\infty) \cap \mathbb{Z}$.

\begin{lemma}
\label{lem:stein:urn}
Let $p$ be the probability mass function of the number $S_n^{\alpha,\beta,m}$ of white balls drawn from the P\'olya-Eggenberger urn by time $n$. Then a random variable $S$ has probability mass function $p$ if and only if for all functions $f \in {\cal{F}}(p)$
\bea \label{Steinpolya}
\e \left[ S( \beta/m + n - S) \Delta f(S-1) + \left\{(n-S) (\alpha/m + S) - S (\beta/m + n - S)\right\}  f(S)  \right] = 0.
\ena
\end{lemma}

We prove Lemma \ref{lem:stein:urn} by applying a general technique for constructing equations such as (\ref{Steinpolya}) from discrete probability mass functions which is of independent interest, see \cite{leyswan}.
We begin with Proposition \ref{steindiscrete} below, a discrete version of the density approach to the Stein equation.

\begin{proposition} \label{steindiscrete}
Let $Z$ have probability mass function $p$ with support
the integer interval $I$, and let
$\psi(k)$ be given by (\ref{def.psi.discrete}) for $k \in I$. If a random variable $X$ with support $I$ has mass function $p$ then for all $f \in {\cal{F}}(p)$,
\bea \label{steindiscreteeq}
\e (\Delta f(X-1 )+\psi (X) f(X))= 0 .
\ena
Conversely, if (\ref{steindiscreteeq}) holds for all $f(k)={\bf 1}(k=\ell), \ell \in I$ then $X$ has mass function $p$.
\end{proposition}

\begin{remark}
\begin{enumerate}
\item The functions $f(k)={\bf 1}(k=\ell), \ell \in I$ are indeed in ${\cal{F}}(p)$, so Proposition \ref{steindiscrete} implies that a random variable $X$ with support $I$ has mass function $p$ if and only if  (\ref{steindiscreteeq}) holds  for all $f \in {\cal{F}}(p)$.
\item
The statement in Proposition \ref{steindiscrete} is  equivalent to Theorem 1.1 given in \cite{leyswan} under a different assumption, namely that equality (\ref{steindiscreteeq}) holds with $g$ replacing $f$ for
all functions
for which $\sum_{k \in I} \Delta(g(k) p(k)) =0$.
 We note that their set-up would translate to test functions $f(k) = g(k+1)$.
\end{enumerate}
\end{remark}

\proof Let $p$ be a real valued function defined on the integer interval $[a,b+1] \cap \mathbb{Z}$ for $a,b \in \mathbb{Z}$ with $a \le b$. Applying the summation by parts formula in the first line below, we obtain
\bea
\label{summation.by.parts}
\sum_{k=a}^b f (k)  \Delta p (k) &=& - \sum_{k=a}^b  p(k+1) \Delta f(k) + f(b+1)p(b+1)-f(a) p(a) \\
&=& - \sum_{k=a+1}^{b+1}  p(k) \Delta f(k-1) + f(b+1)p(b+1)-f(a) p(a) \nonumber \\
&=&  - \sum_{k=a}^b  p(k) \Delta f(k-1) + f(b)p(b+1)-f(a-1) p(a) \label{Abel.right} .
\ena
If  $p(b+1)=0$ and $f(a-1)=0$ then we obtain
\bea
\sum_{k=a}^b f (k)  \Delta p (k)
&=&  - \sum_{k=a}^b p(k) \Delta f(k-1) \label{finite.a.b}.
\ena
Hence
(\ref{steindiscreteeq}) holds when $p$ is a probability mass function with support $[a,b] \cap \mathbb{Z}$ and $f(a-1)=0$.

The case where $I$ is an infinite integer interval follows by applying
Abel's Lemma on summation by parts as modified by \cite{chu2}. In particular, if either \ref{Abel.right}) or the left hand side of (\ref{summation.by.parts})
is convergent upon replacing $b$ by $\infty$, and if
 $\lim_{n \rightarrow \infty} f(n)p(n+1)=0$, then (3a) of \cite{chu2} shows that (\ref{finite.a.b}) holds upon replacing $b$ by infinity,
completing the argument when $I=[a,\infty)$. Similarly, (3b) of \cite{chu2} can be used to argue the case when $I=\mathbb{Z}$.

For $I=(-\infty,b]$, (\ref{Abel.right}) gives that
\beas
\e \psi(X) f(X) &=& \sum_{k=-\infty}^b f(k)\Delta p(k)
= \lim_{a \rightarrow -\infty} \sum_{k=a}^b f (k)\Delta p(k) \\
&=&   \lim_{a \rightarrow -\infty} \left( - \sum_{k=a}^b  p(k) \Delta f(k-1) + f(b)p(b+1)-f(a-1) p(a)\right) .
\enas
Since $p(b+1)=0$
and $p(a)f(a-1) \rightarrow 0$ as $a \rightarrow -\infty$, we obtain that
\beas
\e \psi(X) f(X) &=&  - \lim_{a \rightarrow -\infty}  \sum_{k=a}^b  p(k) \Delta f(k-1) = - \e \Delta f(X-1 )
\enas
and (\ref{steindiscreteeq}) holds.

Conversely, if $X$ with support $I$ satisfies \eqref{steindiscreteeq}
for all functions $f(k) = {\bf 1}(k=\ell)$ for $\ell \in I$, then
\beas
0 &=& \e ( \Delta f(X-1 )+\psi (X) f(X)) \\
&=& \sum_{k \in I} \p(X=k) \left\{ (f(k) - f(k-1) )  + \psi(k) f(k) \right\} \\
&=& \p(X=\ell) - \p(X= \ell + 1)   +\p(X=\ell) \left( \frac{\p(Z=\ell+1)-\p(Z=\ell)}{\p(Z=\ell)} \right),
\enas
and rearranging gives
\beas
\frac{\p( X= \ell + 1)}{\p(X=\ell)} = \frac{\p( Z= \ell + 1)}{\p(Z=\ell)}.
\enas
Hence, if $I$ is $[a,b]$ or $[a,\infty)$ we obtain that for all $j \in I$,
\beas
\frac{\p( X=j)}{\p(X=a)} = \prod_{\ell=a}^{j-1} \frac{\p( X= \ell + 1)}{\p(X=\ell)} = \prod_{\ell=a}^{j-1} \frac{\p( Z= \ell + 1)}{\p(Z=\ell)} = \frac{\p( Z=j)}{\p(Z=a)}.
\enas
Summing over $j \in I$ yields $P(X=a)=P(Z=a)$, whence $P(X=j)=P(Z=j)$ for all $j \in I$. Similarly one may handle the remaining case where $I=(-\infty,b]$.
\bbox

Given a characterization produced by Proposition \ref{steindiscrete}, the following corollary produces varieties of characterizations for the same distribution, each one corresponding to a choice of a function $c$ possessing certain mild properties.

\begin{corollary} \label{notsoverynewdiscretemethod}
Let $Z$ be a random variable with probability mass function $p$ with support $I=[a,b] \cap \mathbb{Z}$ where $a \le b, a, b \in  \mathbb{Z}$, and let $\psi$ be given by (\ref{def.psi.discrete}). Then for all functions $c:[a-1, b] \cap \mathbb{Z} \rightarrow \mathbb{R} \setminus \{0\}$,
a random variable $X$ with support $I$ has mass function $p$ if and only if
\bea \label{need.to.exist}
\e \left[c(X-1) \Delta f(X-1 )+ [ c(X) \psi (X) + c(X) - c(X-1)] f(X) \right]= 0
\ena
for all functions $f \in {\cal{F}}(p)$.
\end{corollary}

\proof Suppose that (\ref{need.to.exist}) holds for $X$ for all functions in $f \in {\cal{F}}(p)$. Since $c(\ell) \not = 0$ for all $\ell \in I$ the functions $f(k) = c(\ell)^{-1} {\bf{1}}(k=\ell), \ell \in I$ all lie in ${\cal{F}}(p)$, so $X$ has mass function $p$ by Proposition \ref{steindiscrete}. For the converse, note that $f \in {\cal{F}}(p)$ implies that $cf \in {\cal{F}}(p)$. Hence if $X$ has mass function $p$ then (\ref{need.to.exist}) follows from Proposition \ref{steindiscrete} upon replacing $f(x)$ by $c(x)f(x)$.
\bbox

Proposition \ref{steindiscrete} applied to $Z$ with mass function $p$ of the Poisson distribution ${\cal P}(\lambda)$, having $\psi(k) = \lambda/(k+1) -1$ by (\ref{def.psi.discrete}) for all $k \in \mathbb{N}_0$,
yields that for all functions $f \in {\cal{F}}(p)$,
\beas  
\e \Delta f(Z-1) = \e \left( 1 - \frac{\lambda}{Z+1} \right)f(Z), \qmq{that is,} \e f(Z-1) = \e \left(\frac{\lambda}{Z+1} \right)f(Z),
\enas
a nonstandard version of a characterization of the Poisson. An extension of Corollary \ref{notsoverynewdiscretemethod} to the case of infinite support produces the usual characterization by the choice $c(k)=k+1$ and the substitution $g(k)=f(k-1)$. Naturally, additional characterizations are produced when using different choices of $c$.

\begin{remark} \label{rem:cond.arb}
When $I=[a,b]\cap \mathbb{Z}$
 we automatically have $\psi(b)=-1$, in which case $c(b)$ does not appear in (\ref{need.to.exist}), and may be assigned a value arbitrarily.
\end{remark}

We now apply Corollary \ref{notsoverynewdiscretemethod} to the distribution of $S_n^{\alpha,\beta,m}$, the number of white balls drawn from the P\'olya-Eggenberger urn by time $n$. We suppress $\alpha,\beta$ and $m$ for notational ease unless clarity demands it. It is well known, and not difficult to verify, that the distribution $p_k=P(S_n=k), k \in \mathbb{Z}$ satisfies
\bea \label{pk.beta.mass}
p_k = {n \choose k} \frac{(\alpha/m)_k (\beta/m)_{n-k}}{(\alpha/m + \beta/m)_n},
\ena
where $(x)_0=1$ and otherwise $(x)_k=x(x+1)\cdots(x+k-1)$ is the rising factorial. The distribution (\ref{pk.beta.mass}) is also known as the beta-binomial and the negative hypergeometric distribution, see \cite{Wilcox}. We now have the ingredients to prove Lemma \ref{lem:stein:urn}.

\medskip
\noindent {\em Proof of Lemma \ref{lem:stein:urn}:}
Taking differences in (\ref{pk.beta.mass}) for $k=0, \ldots, n-1$ yields
\beas
\Delta p_k = {n \choose k} \frac{(\alpha/m)_k (\beta/m)_{n-k}}{(\alpha/m + \beta/m)_n}\left\{\frac{(n-k)(\alpha/m + k) }{ (k+1) (\beta/m + n - k - 1)} -1 \right\},
\enas
while for $k=n$,
\beas
\Delta p_n = - p_n.
\enas
Hence with $\psi(k)=\Delta p_k/p_k$ as in (\ref{def.psi.discrete}) we obtain for $k=0, \ldots, n-1$
\beas
\psi(k) = \frac{(n-k) (\alpha/m + k) - ( k+1) (\beta/m + n - k - 1)}{(k+1) ( \beta/m + n - k - 1)},
\enas
and $\psi(n)=-1$.

In applying Corollary \ref{notsoverynewdiscretemethod}, as $\psi(n)=-1$ we may take the value $c(n)$ arbitrarily, see Remark \ref{rem:cond.arb}. In particular, taking
$c(k) = (k+1) (\beta/m + n - k - 1)$  for all $k=0,\ldots,n-1$ and $c(n) =n$  we obtain (\ref{Steinpolya}). \bbox

\medskip
The next lemma is instrumental in calculating the higher moments of $S_n^{\alpha,\beta,m}$. We let $[x]_0=1$, and otherwise set $[x]_k=x(x-1)\cdots(x-k+1)$, the falling factorial.

\begin{lemma} \label{beta.moment:lem}
For all nonnegative integers $n,a$ and $b$, we have
\bea \label{exp.S.n-S}
\e \left( [S_n^{\alpha,\beta,m}]_a [n-S_n^{\alpha,\beta,m}]_b\right) = \frac{[n]_{a+b}(\alpha/m)_a (\beta/m)_b}{(\alpha/m+\beta/m)_{a+b}}.
\ena
\end{lemma}

\proof First we note that both sides of (\ref{exp.S.n-S}) are zero when $a+b \ge n+1$. This is clear for the right hand side, as the falling factorial $[n]_{a+b}$ is zero. For the left hand side, if $S_n \le a-1$ then $[S_n]_a=0$. On the other hand, if $S_n \ge a$ then $b-1 \ge n-a \ge n-S_n$, in which case $[n-S_n]_b$ is zero.

Now assume $n \ge a+b$. For any $k=0,1,\ldots, n$ we have
\beas
\lefteqn{[k]_a[n-k]_b P(S_n^{\alpha,\beta,m}=k)}\\
&=& [k]_a[n-k]_b {n \choose k}\frac{(\alpha/m)_k(\beta/m)_{n-k}}{(\alpha/m+\beta/m)_n}\\
&=& \frac{[n]_{a+b}(\alpha/m)_a (\beta/m)_b}{(\alpha/m+\beta/m)_{a+b}} {n-a-b\choose k-a}\frac{(\alpha/m+a)_{k-a} (\beta/m+b)_{n-k-b}}{(\alpha/m+\beta/m+a+b)_{n-a-b}}\\
&=& \frac{[n]_{a+b}(\alpha/m)_a (\beta/m)_b}{(\alpha/m+\beta/m)_{a+b}} P(a+S_{n-a-b}^{\alpha+am ,\beta+bm,m}=k).
\enas
Summing over $k=0,1,\ldots,n$ and using that the support of $S_r$ is $\{0,\ldots,r\}$ yields (\ref{exp.S.n-S}). \bbox

\medskip
If $Z$ has the limiting beta distribution ${\cal B}(\alpha/m,\beta/m)$ with density (\ref{Beta.density}), using (\ref{exp.S.n-S}) we obtain
\begin{multline} \label{momentratio}
\e \left( \frac{[S_n]_a [n-S_n]_b}{n^{a + b}}\right)  =
\frac{[n]_{a+b}(\alpha/m)_a (\beta/m)_b}{n^{a+b}(\alpha/m+\beta/m)_{a+b}}\\=
\frac{[n]_{a+b}}{n^{a+b}}\frac{B(\alpha/m+a,\beta/m+b)}{B(\alpha/m,\beta/m)}=
\frac{[n]_{a+b}}{n^{a+b}}\e \left( Z^a(1-Z)^b \right)
\end{multline}
that is, the scaled falling factorial moments of $S_n$ and the power moments of $Z$
differ only by factors of order $1/n$.
This observation can be used to provide a proof of convergence in distribution of $W_n=S_n/n$ to $Z$ by the method of moments, but without a bound on the distributional distance.\\[1ex]

\section{Bounds for the P\'olya-Eggenberger urn model}
Theorem \ref{urntobeta} provides an explicit bound in Wasserstein distance of order $O(1/n)$ between the distribution of $W_n$, the fraction of white balls drawn from the urn by time $n$, and the limiting Beta distribution. For approximating a discrete distribution by a continuous one the Wasserstein distance $d_W$ is a typical distance to use, see for example \cite{gibbssu}. For random variables $X$ and $Y$, this distance is given by
\bea \label{def:Wass}
d_W(X,Y)  = \sup_{h \in {\rm Lip}(1)} \left\vert \e [h(X)]- \e [h(Y)] \right\vert,
\ena
where ${\rm Lip}(1)=\{h:|h(x)-h(y)| \le |x-y|\}$, the class of all Lipschitz-continuous functions with Lipschitz constant less than or equal to 1.
The Wasserstein distance defines a metric on the set of probability measures on $(\RR, {\cal{B}}(\RR))$, the set of real numbers equipped with their Borel $\sigma-$field. On this space, convergence under the Wasserstein distance implies weak convergence. On $([0,1], {\cal{B}}([0,1]))$,
the Wasserstein distance metrizes weak convergence, see \cite{gibbssu}.

\begin{remark} \label{optimal}
The function $h(x) = x(1-x){\bf 1}_{\{x \in [0,1]\}}$ is in ${\rm Lip}(1)$, and applying \eqref{momentratio} with $a=b=1$ we obtain that for all $\alpha \ge 1, \beta \ge 1$ and $m \ge 1$,
\beas
d_W(W_n,Z)   \ge |\e (h (W_n) ) - \e h(Z)| &=& \Bvert \left(
\frac{[n]_{2}}{n^{2}} -1 \right) \e \left( Z(1-Z)\right)\Bvert = \frac{1}{n} \frac{\alpha \beta}{(\alpha+\beta+m)(\alpha+\beta)}.
\enas
Thus the $1/n$ order of the bound in Theorem \ref{urntobeta} cannot be improved.
\end{remark}

\begin{remark} \label{dobler}
Theorem 4.3 of \cite{dobler} provides a bound of order $1/n$ for the Beta approximation to the P\'olya-Eggenberger urn for test functions with bounded first and second derivatives using an exchangeable pair coupling.  The results in \cite{dobler} differ from ours in two significant ways. Firstly, the bound in Theorem 4.3 of \cite{dobler} is expressed in terms of two non-explicit constants $C_1,C_2$ that are defined in Proposition 3.8 of \cite{dobler}. Lemma \ref{lem:allbounds} below provides values of $C_1$.
 The lack of an explicit expression for $C_2$ 
  in \cite{dobler} can be explained by the fact that the solution there is given in terms of ratios of functions which are related to incomplete Beta functions, for which a uniform bound would be difficult.

A more important difference between the present work and \cite{dobler} is that expressing the bound of the latter, presently given in terms of twice differentiable functions, in terms of a bound in a metric, say $d_2$, obtained from twice differentiable functions in the same way that Lipschitz functions yield the Wasserstein metric $d_W$, we have that $d_2 \le d_W$ with equality everywhere not holding. Hence Theorem \ref{urntobeta} implies bounds in the $d_2$ metric, while the reverse does not hold.

\end{remark}

In the following we set our test functions $h$ to be zero outside the unit interval $[0,1]$. For $y > 0$ set
\beas 
\Delta_y f(x)=f\left(x+y\right)-f(x),
\enas
and for a real valued function $g$ on $[0,1]$ we let $||g||=\sup_{w \in [0,1]}|g(w)|$, the supremum norm of $g$.
In the following we recall, with the help of Rademacher's Theorem, 
that a function $h$ is in ${\rm Lip}(1)$ if and only it is absolutely continuous with respect to Lebesgue measure with an almost everywhere derivative bounded in absolute value by 1.

Lemma \ref{betasol:lem} below shows that for all $\{\alpha,\beta\} \subset (0,\infty)$ and functions $h$ for which the expectation ${\cal B}_{\alpha,\beta}h$ exists,
\bea \label{betasol}
f_{\alpha, \beta} (w) = \frac{1}{w^\alpha (1-w)^{\beta}} \int_0^w u^{\alpha-1} (1-u)^{\beta - 1} ( h(u) - {\cal B}_{\alpha,\beta}h) du, \quad w \in [0,1]
\ena
is the unique bounded solution of the Stein equation \eqref{stein:beta}. We continue to omit subscripts when the context makes it clear which parameters are used.
In the proof below we invoke Lemma \ref{lem:allbounds} which yields bounds on the supremum norm of the derivative $f'$ in terms of that same norm on the derivative $h'$ of the test function.

\medskip
{\bf{Proof of Theorem \ref{urntobeta}.}} For $h$ a given function in ${\rm Lip}(1)$, let $f=f_{\alpha/m, \beta/m}$ be the solution of the Stein equation \eqref{stein:beta}  given in \eqref{betasol}.
Replacing $f(z)$ by $f(z/n)$ and dividing by $n$ in \eqref{Steinpolya} results in
\beas
0 &=&
\e S_n \left(\frac{\beta}{nm} + 1 - W_n \right) \Delta_{1/n} f \left(W_n -\frac1n \right) \\
&&+ \e \left\{(n-S_n) \left(\frac{\alpha}{nm} + W_n \right) - S_n  \left(\frac{\beta}{nm} + 1 - W_n\right)  \right\}f(W_n)\\
&=& \e \left[n W_n \left( \frac{\beta}{nm} + 1 - W_n\right) \Delta_{1/n} f \left( W_n -\frac1n \right) + \left\{\frac{\alpha}{m}(1-W_n)  - \frac{\beta}{m} W_n \right\} f(W_n) \right].
\enas

Applying this identity in the Stein equation \eqref{stein:beta}, with $\alpha$ and $\beta$ replaced by $\alpha/m$ and $\beta/m$  respectively, and invoking Lemma \ref{lem:allbounds} below to yield the existence and boundedness of $f'$, we obtain
\bea
\lefteqn{\e h(W_n)- {\cal B}h} \nn \\
&=&  \e \left( W_n (1-W_n ) f'(W_n) + \left[ \frac{\alpha}{m} (1- W_n ) - \frac{\beta}{m} W_n  \right] f(W_n ) \right) \nn \\
&=& \e \left( W_n (1-W_n ) f'(W_n) - n W_n \left( \frac{\beta}{nm} + 1 - W_n\right) \Delta_{1/n} f \left(W_n -\frac1n \right)\right)\nn \\
&=& \e \left( W_n (1-W_n ) f'(W_n) - n W_n \left(  1 - W_n\right) \Delta_{1/n} f \left(W_n -\frac1n \right)\right) + R_1, \label{stein.main.term}
\ena
where, using Lemma \ref{beta.moment:lem} to calculate moments, we obtain
\bea \label{1:thm}
|R_1| = \frac{\beta}{m} \left| \e W_n \Delta_{1/n} f \left(W_n -\frac1n \right) \right| \le \frac{\beta}{nm} \| f' \| \e W_n = \frac{\alpha \beta}{nm (\alpha + \beta)} \| f' \| .
\ena
Writing the difference in (\ref{stein.main.term}) as an integral, we have
\bea
\lefteqn{ \e \left( W_n (1-W_n ) f'(W_n) - nW_n \left(  1 - W_n\right) \Delta_{1/n} f \left(W_n -\frac1n \right)\right) }\nn \\
&=& \e W_n (1-W_n )  \left( f'(W_n) - n \int_{W_n -\frac1n }^{W_n} f'(x) dx \right) \nn \\
&=& n \e  \int_{W_n -\frac1n }^{W_n} W_n(1-W_n)(f'(W_n)  -f'(x)) dx \nn \\
&=&  n \e   \int_{W_n -\frac1n }^{W_n} ( W_n (1-W_n ) f'(W_n)  - x(1-x) f'(x)) dx - R_2, \label{stein.main.term.2}
\ena
where
\beas 
R_2 =  n \e   \int_{W_n-\frac1n}^{W_n} (W_n (1-W_n )-x(1-x)) f'(x) dx.
\enas

To handle $R_2$, using that the solution $f$ of the Stein equation equals $0$ for $x \not\in [0,1]$ to obtain the first inequality,

\bea
|R_2|  &=& \left|  n \e   \int_{W_n -\frac1n}^{W_n}  f'(x) \int_x^{W_n} (1 - 2 y) dy dx \right| \nn \\
&\le& \| f' \|   n \e   \int_{W_n -\frac1n}^{W_n}  \int_x^{W_n} dy dx  \nn \\
&=& \| f' \|   n \e   \int_{W_n-\frac1n}^{W_n} (W_n-x) \, dx  \nn \\
&=& \frac{1}{2n}\| f'\|. \label{2:thm}
\ena

For the first term in (\ref{stein.main.term.2}), substituting using the Stein equation \eqref{stein:beta} with $\alpha$ and $\beta$ replaced by $\alpha/m$ and $\beta/m$, respectively, we obtain
\beas
\lefteqn{ n \e   \int_{W_n -\frac1n }^{W_n} ( W_n (1-W_n ) f'(W_n)  - x(1-x) f'(x)) dx  }\\
&=&  n \e   \int_{W_n -\frac1n }^{W_n} \left\{ h(W_n) - h(x) + \frac{1}{m}\left[(\beta W_n - \alpha( 1-W_n)) f(W_n )   -
(\beta x- \alpha ( 1-x)) f(x) \right] \right\} dx \\
&=&  n \e   \int_{W_n -\frac1n }^{W_n}
\left\{ \int_x^{W_n} h'(y) dy\, + \frac{1}{m}\int_x^{W_n} \left[ \beta y- \alpha ( 1-y) f(y)\right]'dy \right\} dx \\
&=&   n \e   \int_{W_n -\frac1n }^{W_n} \left\{ \int_x^{W_n} h'(y) dy\,  + \frac{1}{m}\int_x^{W_n}
\left[\beta y - \alpha ( 1-y)) f'(y ) + (\beta + \alpha)f(y)  \right] dy  \right\} dx .
\enas

We bound the inner integrals separately. Firstly,
\bea \label{3:thm}
\left|  n \e   \int_{W_n-\frac1n}^{W_n}  \int_x^{W_n} h'(y) dy dx \right|
&\le& n|| h' ||  \e   \int_{W_n-\frac1n}^{W_n}   (W_n - x)  dx = \frac{1}{2n} ||h'||.
\ena

Next, recalling that $ 0 \le W_n \le 1$ and noting that $|(\beta y - \alpha ( 1-y))| \le \alpha \vee \beta$ for $0 \le y \le 1$,
\bea
\lefteqn{ \left|  \frac{n}{m} \e   \int_{W_n-\frac1n}^{W_n} \int_x^{W_n} (\beta y - \alpha ( 1-y)) f'(y ) dy \, dx\right|} \nn \\
&=&  \left| \frac{n}{m} \e   \int_{W_n-\frac1n}^{W_n}  (\beta y - \alpha ( 1-y)) f'(y )  \int_{W_n-\frac1n}^{y} \, dx \, dy \right|\nn \\
&\le& \frac{n}{m} || f' || (\alpha \vee \beta) \e   \int_{W_n-\frac{1}{n}}^{W_n} \int_{W_n-\frac{1}{n}}^{y} \, dx \, dy \nn \\
&=& \frac{1}{2nm} \| f' \| ( \alpha \vee \beta ). \label{4:thm}
\ena
Arguing in a similar fashion, we obtain
\bea
\left|\frac{n}{m} \e   \int_{W_n -\frac1n }^{W_n}  \int_x^{W_n}  (\beta + \alpha ) f(y)  dy  dx  \right|
\le \frac{n}{m}\|f\|(\beta + \alpha ) \e   \int_{W_n-\frac1n}^{W_n}  \int_{W_n-\frac1n}^y dx dy =
\frac{1}{2nm}  \|f\| (\alpha + \beta). \label{5:thm}
\ena

Collecting the bounds (\ref{1:thm}), (\ref{2:thm}), (\ref{3:thm}), (\ref{4:thm})  and (\ref{5:thm}) yields
\beas 
\left| \e h(W_n) - \e h(W) \right| & \le& \left( \frac{  m + \alpha \vee \beta}{2nm} + \frac{\alpha \beta}{nm(\alpha+\beta)} \right) \| f' \|  + \frac{1}{2n} ||h'||
 + \frac{1}{2nm}  (\alpha + \beta) \| f \|.
\enas
Note that by  Lemma \ref{lem:dobler}
\beas
\frac{1}{2nm}  (\alpha + \beta) \| f \| \le \frac{1}{2nm} \left( \frac{2}{\alpha/m + \beta/m} \right)  ||h'|| = \frac1n ||h'||.
\enas
The theorem now follows by invoking  Lemma \ref{lem:allbounds}.
\bbox

\begin{lemma}
\label{betasol:lem}
For any $\{\alpha,\beta\} \subset (0,\infty)$ and real valued function $h$ on $[0,1]$ such that the expectation ${\cal B}_{\alpha, \beta} h $
of $h$ exists, the function $f$ given by (\ref{betasol})
is the unique bounded solution of (\ref{stein:beta}).
\end{lemma}

\proof It is straightforward to verify that $f$ as given in (\ref{betasol}) is a solution of (\ref{stein:beta}).
Writing the associated homogeneous equation as
\beas
(w^{\alpha-1}(1-w)^{\beta-1})^{-1}(w^\alpha (1-w)^\beta g(w))'=0
\enas
we find that all solutions to (\ref{stein:beta}) are given by
\beas 
f(w) + cg(w), \quad \mbox{for some $c \in \mathbb{R}$, where} \,\,  g(w)=\frac{1}{w^\alpha(1-w)^\beta}.
\enas
The claim follows since $g(w)$ is unbounded at the endpoints of the unit interval for all $c \not =0$, and Lemma \ref{lem:allbounds} below demonstrates that $f(w)$ is bounded. \bbox

\medskip
Since the expectation of $h(Z)-{\cal B}h$ is zero when $Z \sim {\cal B}(\alpha,\beta)$, we may also write
\bea \label{alternativesol}
f(w) = -  \frac{1}{w^\alpha (1-w)^{\beta}} \int_w^1 u^{\alpha-1} (1-u)^{\beta - 1} ( h(u) - {\cal B}h) du.
\ena

From Proposition 3.8 in \cite{dobler} we quote the following result.

\begin{lemma}\label{lem:dobler}
The solution $f$, given in (\ref{betasol}),  of (\ref{stein:beta})  for $h$ a Lipschitz function on $[0,1]$ satisfies
\beas 
|| f|| \le \frac{2}{\alpha + \beta}  || h'||.
\enas
\end{lemma}

The cases in the bounds of Lemma \ref{lem:allbounds} reflect the behaviour of the function
\bea \label{def:g}
g(w)=w^{\alpha-1}(1-w)^{\beta-1} \qmq{for $w \in [0,1]$,}
\ena
as described in Lemma \ref{monotone}. In the following we will use the terms decreasing and increasing in the non-strict manner, for example, a constant function is both increasing and decreasing.
Let
\beas
x_{\alpha,\beta}=\frac{\alpha-1}{\alpha+\beta-2}.
\enas

\begin{lemma}\label{monotone}
For $\{\alpha,\beta\} \subset (-1,\infty)$, the function $g:[0,1] \rightarrow [0,\infty)$ given in (\ref{def:g})
has the following behaviour.

\msk
\begin{tabular}{|c|c|c|c| } \hline
& $\beta < 1$ &   $\beta = 1$ &  $\beta > 1$ \\ \hline
 $\alpha< 1$ & decreasing on $[0, x_{\alpha,\beta}]$ & decreasing & decreasing \\
 & increasing on $[x_{\alpha,\beta}, 1]$ & & \\ \hline
 $\alpha = 1 $ & increasing & constant & decreasing \\ \hline
 $\alpha > 1$ & increasing & increasing & increasing on $[0, x_{\alpha,\beta}]$ \\
 & & & decreasing on $[x_{\alpha,\beta}, 1]$ \\  \hline
\end{tabular}

\end{lemma}

\msk
\proof Clearly when $\alpha=1$ and $\beta=1$ the function $g(w)$ is constant. Otherwise,
taking derivative in (\ref{def:g}) yields
\beas
g'(w) &=& (\alpha -1) w^{\alpha -2} (1-w)^{\beta-1}- (\beta-1) w^{\alpha-1} (1-w)^{\beta-2}\\
&=&  w^{\alpha -2} (1-u)^{\beta-2} \{ (\alpha -1) (1-w) - (\beta - 1) w \} .
\enas
The expression is non-negative if and only if
\bea \label{g'ge0}
(\alpha -1) (1-w) \ge (\beta - 1) w.
\ena
When $\alpha \ge 1$ and $\beta \le 1$ inequality (\ref{g'ge0}) is always satisfied. Similarly (\ref{g'ge0}) holds with the non-strict inequality reversed when $\alpha \le 1$ and $\beta \ge 1$.
The remaining two cases $\alpha < 1, \beta < 1$ and $\alpha > 1, \beta > 1$ follow by solving the inequality.
\bbox

Our next result bounds the magnitude of the derivative of the solution $f$ in terms of $h$.

\begin{lemma}
\label{lem:allbounds} For $\{\alpha,\beta\} \subset (0,\infty)$ let $f=f_{\alpha,\beta,h}$ be the solution to \eqref{stein:beta} given by \eqref{betasol} for an
absolutely continuous function $h$. Then
\bea \label{f.bound.when.h.Lip}
||f'||
\le    b_0||h-{\cal B}h|| +  b_1||h'|| \le  (b_0 + b_1)||h'||,
\ena
where $b_0=b_0(\alpha, \beta)$ and $b_1=b_1(\alpha, \beta)$
are given by
\beas
b_0 = \left\{ \begin{array}{ll}
4 \max  \left(|\alpha -1|;  |\beta -1| \right) & \mbox{ if }   \alpha \le 2, \beta \le  2;\\
(\alpha + \beta - 2)^2 \max  \left( \frac{\alpha-1}{(\alpha -2)^2}; \frac{|\beta -1|}{\beta^2}  \right)
& \mbox{ if }  \alpha > 2, \beta \le 2;\\
 (\alpha + \beta -2)^2 \max \left(\frac{|\alpha-1|}{\alpha^2}; \frac{\beta-1}{(\beta-2)^2} \right)
  & \mbox{ if }   \alpha \le 2, \beta >2  ;\\
  (\alpha + \beta - 2)^2 \max\left(\frac{1}{\alpha-1}; \frac{1}{\beta-1} \right)
 & \mbox{ if }  \alpha > 2, \beta > 2\\
\end{array}
\right.
\enas
and
\beas
b_1 = \left\{ \begin{array}{ll}
4 \left( 1 + \frac{\max  \left(\alpha ;  \beta \right)}{\alpha + \beta} \right) & \mbox{ if }   \alpha \le 2, \beta \le  2;\\
\frac{(\alpha + \beta -2)^2 }{\min( \alpha-2, \beta )^2}
+ 2  \max \left( \frac{\alpha}{\alpha-2}; 1  \right)
 & \mbox{ if } \alpha > 2, \beta \le 2;\\
 \frac{(\alpha + \beta -2)^2 }{\min( \alpha, \beta-2)^2} +2 \max \left( 1; \frac{\beta}{\beta-2} \right)
 & \mbox{ if }   \alpha \le 2, \beta >2  ;\\
 \frac{(\alpha + \beta -2)^2 }{\min( \alpha-1, \beta-1 )^2} + 2  \max \left(\frac{\alpha}{\alpha-1}; \frac{\beta}{\beta-1} \right)
 & \mbox{ if }  \alpha > 2, \beta > 2.\\
\end{array}
\right.
\enas

\end{lemma}
\bigskip

\proof By replacing $h$ by $h-{\cal B}_{\alpha,\beta}h$ we may assume ${\cal B}_{\alpha,\beta}h=0$.
 Rewriting the Stein equation \eqref{stein:beta} yields
\beas
w(1-w)f'(w) = h(w) +(\beta w -\alpha(1-w))f(w),
\enas
so to show  (\ref{f.bound.when.h.Lip}) it suffices to demonstrate that for all $w \in [0,1]$
\beas
|h(w) +(\beta w -\alpha(1-w))f(w)| \le (b_0||h||+b_1||h'||) w (1-w).
\enas

Using (\ref{betasol}) and integration by parts we obtain
\begin{multline*}
\alpha(1-w) f(w)
= h(w) + \frac{\beta-1}{w^{\alpha} (1-w)^{\beta-1}} \int_0^w  u^{\alpha} (1-u)^{\beta-2} h(u)  du \\
-\frac{1}{w^{\alpha} (1-w)^{\beta-1}} \int_0^w u^{\alpha} (1-u)^{\beta-1} h'(u)  du
\end{multline*}
and, now applying (\ref{alternativesol}),
\begin{multline*}
\beta w f(w)
= - h (w) - \frac{\alpha-1}{w^{\alpha-1} (1-w)^{\beta}} \int_w^{1}  u^{\alpha-2} (1-u)^{\beta} h(u)  du  \\
- \frac{1}{w^{\alpha-1} (1-w)^{\beta}}  \int_w^{1}  u^{\alpha-1} (1-u)^{\beta} h'(u)  du.
\end{multline*}
Hence
\begin{multline*}
h(w) +(\beta w -\alpha(1-w))f(w)\\
= -\frac{\beta-1}{w^{\alpha} (1-w)^{\beta-1}} \int_0^w  u^{\alpha} (1-u)^{\beta-2} h(u)  du
+\frac{1}{w^{\alpha} (1-w)^{\beta-1}} \int_0^w u^{\alpha} (1-u)^{\beta-1} h'(u)  du
+ \beta w f(w) \\
= - \frac{\alpha-1}{w^{\alpha-1} (1-w)^{\beta}} \int_w^{1}  u^{\alpha-2} (1-u)^{\beta} h(u)  du
 - \frac{1}{w^{\alpha-1} (1-w)^{\beta}}  \int_w^{1}  u^{\alpha-1} (1-u)^{\beta} h'(u)  du
- \alpha(1-w) f(w).
\end{multline*}

From Lemma \ref{lem:dobler} we immediately have the bounds
\bea  \label{imm1}
|\beta w f(w) | \le \frac{2\beta}{\alpha + \beta} w || h' || \qmq{and} | \alpha(1-w) f(w) | \le \frac{2 \alpha}{\alpha + \beta}   (1-w) ||h'||.
\ena

As $0 \le u \le 1$,
\begin{multline*}
\left|- \frac{\beta-1}{w^{\alpha} (1-w)^{\beta-1}} \int_0^w  u^{\alpha} (1-u)^{\beta-2} h(u)  du  +\frac{1}{w^{\alpha} (1-w)^{\beta-1}} \int_0^w u^{\alpha} (1-u)^{\beta-1} h'(u)  du \right|  \\
\le   \frac{|| h ||\,|\beta-1|}{w^{\alpha} (1-w)^{\beta-1}} \int_0^w  u^{\alpha} (1-u)^{\beta-2}  du
 +  \frac{|| h'||}{w^{\alpha} (1-w)^{\beta-1}} \int_0^w u^{\alpha} (1-u)^{\beta-1}   du\\
\le (  | \beta -1| || h|| + || h'|| )  \frac{1}{w^{\alpha} (1-w)^{\beta-1}} \int_0^w  u^{\alpha} (1-u)^{\beta-2}  du.
\end{multline*}

When $u^{\alpha} (1-u)^{\beta-2} $ is increasing on $[0,x_*]$ then for $w \in [0,x_*]$ we can bound this expression by
\beas
|| h || w (1-w)   \frac{|\beta -1|}{(1-x_*)^2} + || h'|| w (1-w)\frac{1}{(1-x_*)^2},
\enas
and now using the first inequality in (\ref{imm1}), we obtain
\beas
\lefteqn{\left|  h(w) +(\beta w -\alpha(1-w))f(w) \right|}\\
&\le&  || h || w (1-w)  \frac{|\beta -1|}{(1-x_*)^2}+ || h'|| w (1-w) \left( \frac{1}{(1-x_*)^2} +   \frac{2 \beta}{(\alpha + \beta)(1-x_*)}\right)\\
&\le& \frac{|\beta -1|}{(1-x_*)^2} || h || w (1-w)  +  \left( \frac{1}{(1-x_*)^2}+ \frac{2 \beta }{(\alpha + \beta)(1-x_*)} \right) || h'|| w (1-w).
\enas

Similarly,
\begin{multline*}
\left|- \frac{\alpha-1}{w^{\alpha-1} (1-w)^{\beta}} \int_w^{1}  u^{\alpha-2} (1-u)^{\beta} h(u)  du
- \frac{1}{w^{\alpha-1} (1-w)^{\beta}}  \int_w^{1}  u^{\alpha-1} (1-u)^{\beta} h'(u)  du \right|\\
\le ( | \alpha-1| || h|| + || h'|| )  \frac{1}{w^{\alpha-1
} (1-w)^{\beta}} \int_w^{1}  u^{\alpha-2} (1-u)^{\beta}  du
\end{multline*}
and if $u^{\alpha-2} (1-u)^{\beta} $ is decreasing on $[x_*,1]$ then for $w \in [x_*,1]$ we can bound this expression by
\beas
|| h || w (1-w) \frac{|\alpha -1|}{x_*^2} + || h'|| w (1-w)  \frac{1}{x_*^2}.
\enas
Now using the second inequality in (\ref{imm1}), we obtain
\beas
\lefteqn{\left|  h(w) +(\beta w -\alpha(1-w))f(w) \right|}\\
&\le& || h || w (1-w)   \frac{|\alpha -1|}{x_*^2}   + \left(\frac{1}{x_*^2} +  \frac{2 \alpha}{(\alpha + \beta) x_*}  \right) || h'|| w (1-w).
\enas
Hence we may take
\bea \label{b0:general}
b_0=\max\left( \frac{|\alpha -1|}{x_*^2}, \frac{|\beta -1|}{(1-x_*)^2}\right)
\ena
and
\bea \label{b1:general}
b_1=\max\left\{\frac{1}{x_*^2} +  \frac{2 \alpha}{(\alpha + \beta) x_*}  ; \frac{1}{(1-x_*)^2}+\frac{2 \beta}{(\alpha + \beta)(1-x_*)}
\right\}
\ena

In view of Lemma \ref{monotone} we  distinguish four cases.

\medskip
{\bf{Case 1}}. $\alpha \le 2, \beta \le 2$. By Lemma \ref{monotone}, $u^{\alpha} (1-u)^{\beta-2}$ is increasing and $  u^{\alpha-2} (1-u)^{\beta}  $ is decreasing. Setting $x_*=1/2$, by (\ref{b0:general}) and (\ref{b1:general}) we obtain
\beas
b_0 = 4 \max  \left( |\alpha -1|; |\beta -1| \right) \qmq{and}
b_1=4 + \frac{4}{\alpha + \beta}  \max  \left(\alpha ;  \beta \right) = 4 \left( 1 + \frac{\max  \left(\alpha ;  \beta \right)}{\alpha + \beta} \right).
\enas

\msk
{\bf{Case 2}}. $\alpha > 2, \beta \le 2$. In this case, from Lemma \ref{monotone}, $u^{\alpha} (1-u)^{\beta-2}$ is increasing, and $u^{\alpha-2} (1-u)^{\beta}$ is decreasing on $[x_{\alpha-1,\beta+1},1]$. Setting $x_*=x_{\alpha-1,\beta+1}$ and noting that $x_{\alpha,\beta}+x_{\beta,\alpha}=1$, by (\ref{b0:general}) and (\ref{b1:general}) we obtain
$$ b_0 = \max  \left(\frac{\alpha-1}{x_{\alpha-1,\beta+1}^2} ;  \frac{|\beta -1|}{x_{\beta+1,\alpha-1}^2}
 \right) = (\alpha + \beta - 2)^2 \max  \left( \frac{\alpha-1}{(\alpha -2)^2}; \frac{|\beta -1|}{\beta^2}  \right) $$
and
\beas
b_1&= &\max\left\{\frac{1}{x_{\alpha-1,\beta+1}^2} +  \frac{2 \alpha}{(\alpha + \beta) x_{\alpha-1,\beta+1}}  ; \frac{1}{x_{\beta+1,\alpha-1}^2}+\frac{2 \beta}{(\alpha + \beta) x_{\beta+1,\alpha-1}}
\right\}    \\
&\le&
\frac{(\alpha + \beta -2)^2 }{\min( \alpha-2, \beta )^2}+ 2 \left(\frac{\alpha + \beta -2}{\alpha + \beta}\right)  \max \left( \frac{\alpha}{\alpha-2}; 1  \right),
\enas
and bounding $(\alpha + \beta -2)/(\alpha + \beta) $ by 1 gives the assertion.

\msk
{\bf{Case 3}}. $\alpha \le 2, \beta >  2$.
In this case, from Lemma \ref{monotone}, $ u^{\alpha} (1-u)^{\beta-2} $ is increasing on $[0, x_{\alpha+1, \beta-1}]$, and  $  u^{\alpha-2} (1-u)^{\beta}  $ is decreasing. Setting $x_*=x_{\alpha+1, \beta-1}$, by (\ref{b0:general}) and (\ref{b1:general}) we obtain
$$ b_0 = \max  \left( \frac{|\alpha -1|}{x_{\alpha+1,\beta-1}^2}; \frac{\beta -1}{x_{\beta-1,\alpha+1}^2}
 \right) = (\alpha + \beta -2)^2 \max \left(\frac{|\alpha-1|}{\alpha^2}; \frac{\beta-1}{(\beta-2)^2} \right) $$
and
\beas
b_1& =&  \max \left\{\frac{1}{x_{\beta-1,\alpha+1}^2}+ \frac{2\beta}{(\alpha + \beta)x_{\beta-1,\alpha+1}} ; \frac{1}{x_{\alpha+1,\beta-1}^2} + \frac{2 \alpha}{(\alpha + \beta)x_{\alpha+1,\beta-1}} \right\} \\
&\le&
\frac{(\alpha + \beta -2)^2 }{\min( \alpha, \beta-2)^2} +2 \left(\frac{\alpha + \beta -2}{\alpha + \beta}\right) \max \left( 1; \frac{\beta}{\beta-2} \right)  .
\enas

\msk
{\bf{Case 4}}. $\alpha > 2, \beta > 2$.
In this case, from Lemma \ref{monotone}, $ u^{\alpha} (1-u)^{\beta-2} $ is increasing on $[0, x_{\alpha+1, \beta-1}]$, and  $  u^{\alpha-2} (1-u)^{\beta}  $ is decreasing on $[x_{\alpha-1, \beta+1}, 1]$.
Noting that
$$
 x_{\alpha-1,\beta+1} < x_{\alpha, \beta} <  x_{\alpha+1,\beta-1},
$$
setting $x_*=x_{\alpha, \beta}$, by (\ref{b0:general}) and (\ref{b1:general}) we obtain
$$ b_0 = \max  \left(  \frac{\alpha -1}{x_{\alpha,\beta}^2}; \frac{\beta -1}{x_{\beta,\alpha}^2} =
(\alpha + \beta - 2)^2 \max\left(\frac{1}{\alpha-1}; \frac{1}{\beta-1} \right)
 \right)$$
and
\beas
b_1 &= & \max\left\{ \frac{1}{x_{\alpha,\beta}^2} + \frac{2 \alpha}{(\alpha + \beta) x_{\alpha,\beta}} ;  \frac{1}{x_{\beta,\alpha}^2} + \frac{2 \beta}{(\alpha + \beta) x_{\beta,\alpha}}\right\} \\
&\le&
\frac{(\alpha + \beta -2)^2 }{\min( \alpha-1, \beta-1 )^2} + 2 \left(\frac{\alpha + \beta -2}{\alpha + \beta}\right)  \max \left(\frac{\alpha}{\alpha-1}; \frac{\beta}{\beta-1} \right).
\enas

For the final inequality in (\ref{f.bound.when.h.Lip}), with $p(y; \alpha, \beta)$ denoting the ${\cal B}(\alpha,\beta)$ density in (\ref{Beta.density}), we have
\beas
||h -{\cal B}h || &=& \sup_{x \in [0,1]} |h(x)-{\cal B}h| \le \sup_{x \in [0,1]}\int_0^1 |  h(x)-h(y)| p(y; \alpha, \beta) dy  \\
&\le& ||h'|| \sup_{x \in [0,1]} \int_0^1 |x-y|  p(y; \alpha, \beta) dy \le ||h'||. 
\enas

\bbox


\section{Distance of the distribution of $L_{2n}$ to the Arcsine law}
\label{arc:sec}
We rely on \cite{doblerArcSine} for the following argument, noting that in \cite{doblerArcSine} no explicit bound is obtained.

{\it Proof of Theorem \ref{arc:thm}.} Let $p$ be the mass function of $L_{2n}$ given by (\ref{Ln:dist}), and let $Z$ have the Arcsine distribution. Applying Proposition \ref{steindiscrete} for $p$, followed by Corollary \ref{notsoverynewdiscretemethod} with the choice
\beas
\psi(k)=\frac{2k-n+1}{(k+1)(2(n-k)-1)} \qmq{and} c(k)=(k+1)(2(n-k)-1),
\enas
\cite{doblerArcSine} arrives at the version of Lemma \ref{lem:stein:urn}, showing that $W_n=(2n)^{-1}L_{2n}$ is so distributed if and only if
\bea \label{arc.char}
\e\left[nW_n\left(1-W_n+\frac{1}{2n}\right)\Delta_{1/n} f\left(W_n-\frac{1}{n}\right)+\left(\frac{1}{2}-W_n\right)f(W_n)\right]=0
\ena
for all functions $f \in {\cal F}(p)$.

Now following steps as those in Theorem \ref{urntobeta} for the P\'olya urn, collecting the estimates from the proof of Theorem 3.1 of \cite{doblerArcSine}  shows that if $f$ is the solution (\ref{betasol}) to (\ref{stein:beta}) with $\alpha=\beta=1/2$, so that $b_0=2$ and $b_1=6$,  then for all differentiable functions $h$ one has
\beas
\left| \e h(W_n)-\e h(Z) \right| \le \frac{1}{2n}||h'||+ \frac{1}{2n}||f||+ \left( \frac{1}{4n}+\frac{n+2}{2n^2}+\frac{3}{4n}\right)||f'||.
\enas
Applying the bounds of Lemma \ref{lem:allbounds} as well as Lemma \ref{lem:dobler} yields the bound in Theorem \ref{arc:thm}.

\medskip
Applying (\ref{arc.char}) with $f(w)$ replaced by $g(w)=1$ and $g(w)=w$ yields
\beas
\e W_n=\frac{1}{2} \qmq{and} \e W_n^2=\frac{3}{8}+\frac{1}{8n}
\enas
respectively.\ignore{
\beas
\Delta_{1/n}f(w)=(w+\frac{1}{n})-w=\frac{1}{n}
\qmq{hence}
\Delta_{1/n}f(w-\frac{1}{n})=\frac{1}{n}
\enas
so
\begin{multline*}
0=\e\left[ W_n\left(1-W_n+\frac{1}{2n}\right)+\left(\frac{1}{2}-W_n\right)W_n\right]\\
=\e[-2W_n^2+\left(\frac{3}{2}+\frac{1}{2n}\right)W_n]=-2\e[W_n^2]+\frac{3}{4}+\frac{1}{4n}.
\end{multline*}
}
For $f(w)=w^2/2$, a function in ${\rm Lip}(1)$ on $[0,1]$, since $\e Z^2=3/8$ we obtain
\beas
\e f(W_n)-\e f(Z) = \frac{1}{16n}.
\enas
Hence the $O(1/n)$ rate cannot be improved. \bbox

\bigskip
{\bf{Acknowledgements.}} We would like to thank the Keble Advanced Studies Centre, Oxford, for support, and an anonymous referee for helpful comments which lead to an improvement of the paper.

\end{document}